\documentclass[12pt,]{article}
\usepackage[russian]{babel}
\newtheorem{theorem}{\indent{\sc Теорема}}
\newtheorem{lemm}{\indent{\sc Лемма}}
\newtheorem{corollary}{\indent{\sc Следствие}}

\newtheorem{definition}{\indent{\sc определение}}
\begin{document}
\textit{}

УДК 511,517
\textit{}

\textit{}

\begin{center}
\textbf{\sc Гипотеза Римана }

\textbf{ \sc И. Ш. Джаббаров }

\sc Гянджинский университет, Гянджа, Азербайджан

\textit{jabbarovish@rambler.ru} 

\begin{abstract}

\textbf{Аннотация: }В статье дается доказательство знаменитой гипотезы Римана. Метод доказательства основан на приближении дзета-функции на правой половине критической полосы отрезками произведений типа эйлеровских.
\end{abstract}
\end{center}
\vspace{-3mm}

Появление дзета-функции и аналитических методов в теории чисел связано с именем Л. Эйлера (см. [18, стр.54]). В 1748 г. Эйлер ввел дзета-функцию
\begin{equation} \label{1}
\zeta (s)=\sum _{n=1}^{\infty } n^{-s} ,s>1,
\end{equation}
рассматривая ее как функцию вещественного переменного s. Используя тождество

\[\zeta (s)=\prod _{p} \left(1-p^{-s} \right)^{-1} ,s>1,\]
где произведение берется по всем простым, он дал аналитическое доказательство теоремы Евклида о бесконечности множества простых чисел. Эйлер дал соотношение, современная формулировка которого эквивалентна функциональному уравнению Римана (см. [13]).

В 1798 г. Лежандр сформулировал для количества $\pi (x)$, простых, не превосходящих $x$ следующее равенство $\mathop{\lim }\limits_{x\to \infty } \frac{\pi (x)\ln x}{x} =1$ и предположил, что выполняется более точное соотношение: $\pi (x)=x/(\ln x-B(x))$, где при $x\to \infty $ $B(x)$ стремится к постоянному $B=1.083...$.

Еще до Лежандра Гаусс предположил, что $\pi (x)$ приближается, с меньшей ошибкой при помощи функции  $\int _{2}^{x}\frac{du}{\ln u}  $. Согласно этому предположению для $B$ в формуле Лежандра получается значение $B=1$.

В 1837 г. используя и развивая идеи Эйлера, Л. Дирихле дал обобщение теоремы Евклида для арифметических прогрессий, рассматривая L -функции. Дирихле пытался доказать формулу Лежандра вводя понятие асимптотического закона.

В 1851 и 1852 гг. Чебышев получил точные результаты. Он показал, что если при $x\to \infty $ отношение $\pi (x)\ln x/x$ стремится к пределу, то этим пределом будет 1, как и предпологал Лежандр. Он, также, установил, что для постоянной $B$ верным значением может быть только 1. В работах Чебышева исследование функции Эйлера $\zeta (s)$ поднят на более высокий уровень.

Большое значение дзета-функции для аналитической теории чисел был открыт Риманом в 1859-м году. Возможно [41], Риман занялся исследованием дзета-функции под влиянием достижений Чебышева, в своем знаменитом мемуаре [9] он впервые рассмотрел дзета-функцию как функцию комплексного переменного и связал проблему распределения простых чисел с расположением комплексных нулей дзета-функции. Он доказал функциональное уравнение

\[\xi (s)=\xi (1-s);\xi (s)=\frac{1}{2} s(s-1)\pi ^{-s/2} \Gamma (\frac{s}{2} )\zeta (s)\]
и сформулировал несколько гипотез о дзета-функции. Одной из них (далее ГР) суждено было стать центральной проблемой для всей математики. Эта гипотеза утверждает, что все комплексные нули дзета- функции, расположенные в критической полосе $0<Res<1$, лежат на критической прямой $Res=0.5$.

Д. Гильберт в своем докладе в Международном Парижском Конгрессе 1900 г. включил эту гипотеза в список своих 23 математических проблем. Несмотря на попытки математиков, доказать эту гипотезу не удавалось. Чтобы достичь прогресса в направлении доказательств ГР были развиты следующие разделы аналитической теории чисел:

1. Исследовании областей, свободных от нулей дзета-функции;

2. Оценки плотности распределения нулей в критической полосе и их приложения;

3. Изучение нулей на критической прямой;

4. Изучение распределения значений дзета-функции в критической полосе;

5. Вычислительные задачи, связанные с нулями и пр.

Эти направления -- классические и в литературе можно найти достаточно полное изложение как исторических, так и других аспектов вопросов, связанных с ГР (см. [3, 6, 12, 16, 17, 18, 22, 24,40-44]). Здесь, во введении мы затронем вкратце работы 4-го направления и некоторые современные идеи, связанные с ГР.

Изучение распределения значений дзета-функции было начато Г. Бором (см. [25, стр.279]). В работе [2], совместно с Курантом Р., было доказана теорема о всюду плотности значений $\zeta (\sigma +it),-\infty <t<\infty ,\sigma \in (1/2,1)$.

Результаты С.М. Воронина [26-33], связанные со свойством универсальности дзета-функции, подняли на новый уровень исследования дзета-функции Римана и других функций, определяемых рядами Дирихле. В работах С.М. Воронина изучены распределения значений некоторых рядов Дирихле, дано в более общей форме новое решение проблемы Д. Гильберта о дифференциальной независимости дзета-функции и L --функций. О других обобщениях и улучшениях см. ([1, 14-16]).

В последние несколько лет были начаты изучения некоторых семейств рядов Дирихле, целью которых была рассмотрение вопросов распределения нулей дзета- и $L$--функций (см.[40-42]). Б. Багчи рассмотрел (см. [15 - 16]) семейство рядов Дирихле, определяемых, при помощи следующего произведения по всем простым

\[F(s,\theta )=\prod _{p} \left(1-\chi _{p} (\theta )p^{-s} \right)^{-1} ,\]
где $Res>1$, $\theta $ принимает значения из топологического произведения окружностей $|z_{p} |=1,z_{p} \in C$, и $\chi _{p} (\theta )$ является проекцией $\theta $ на окружность $|z_{p} |=1$. Он показал, что эта функция может быть аналитически продолжена в полуплоскость \textit{$Res>1/2$} и не имеет там нулей почти для всех $\theta $. Здесь мера является мерой Хаара. В работах [1,14-16] рассмотрены вопросы, связанные со свойством совместной универсальности некоторых рядов Дирихле. С использованием эргодических методов построены также вероятностные меры.

В работе [11] был дан эквивалентный вариант упомянутого выше результата Б. Багчи рассмотрением функции
\begin{equation} \label{2}
F(s,\theta )=\prod _{p} \left(1-e^{2\pi i\theta _{p} } p^{-s} \right)^{-1} ,0\le \theta _{p} \le 1
\end{equation}
в кубе $\Omega =[0,1]\times [0,1]\times \cdot \cdot \cdot $ с произведением мер Лебега.

 В работах [33 - 38] были изучены вопросы о расстояниях последовательных нулей дзета-функции, расположенных на критической прямой, о числах нулей в кругах небольшого радиуса у близкой окрестности критической прямой, а также, о кратных нулях дзета-функции.

 В настоящей работе мы изучаем распределение специальных кривых вида $(\{ t\lambda _{n} \} )_{n\ge 1} $ (знак $\{ \} $ означает дробную часть и $\lambda _{n} >0,\lambda _{n} \to \infty $, когда $n\to \infty $) в подмножествах бесконечномерного единичного куба, на которых некоторые ряды расходятся. В работах [44- 46] установлено, что в бесконечномерном единичном кубе $\Omega $существует мера, отличная от меры Хаара, не являющаяся инвариантной относительно сдвигов (mod 1) и вышеуказанная кривая неизмерима в смысле этой меры, являясь множеством меры нуль в смысле Хаара. В качестве приложения мы доказываем справедливость Гипотезы Римана. Для установления последней сначала мы будем приближать $\zeta (s)$ в некотором круге, расположенном на правой половине критической полосы с помощью частичных произведений вида (\ref{2}), используя лемму С.М. Воронина (см. лемму 2). Далее, мы будем распространять полученное соотношение на всю правую половину критической полосы, используя специальную структуру множества расходимости некоторых рядов (см. п.5).
\begin{definition}
\textit{. Пусть $\sigma :N\to N$ является произвольным взаимно-однозначным отображением множества натуральных чисел на себя. Если найдется натуральное число $m$ такое, что $\sigma (n)=n$ для любого $n>m$,  тогда мы будем говорить, что $\sigma $ является конечной перестановкой. Подмножество $A\subset \Omega $ мы будем называть конечно-симметричным, если для любого элемента $\theta =(\theta _{n} )\in A$ и любой конечной перестановки $\sigma $ имеем $\sigma \theta =(\theta _{\sigma (n)} )\in A$.}
\end{definition}

Пусть $\Sigma $ обозначает множество всех конечных перестановок. Определим на этом множестве произведение двух конечных перестановок как композиции отображений. Тогда $\Sigma $ становится группой, которая содержит каждую группу подстановок степени $n$ в качестве подгруппы (мы рассматриваем каждую подстановку $\sigma $ степени $n$ как конечную перестановку в смысле определения 1, для которой $\sigma (m)=m$ когда $m>n$). Множество $\Sigma $ является счетным множеством и мы можем расположить его элементы в последовательность.
\begin{theorem}
\textit{. Пусть$0<r<1/4$-действительное число. Тогда существует последователь-ность $(\bar{\theta }_{n} )_{n\ge 1} $ элементов $\Omega \, \, (\bar{\theta }_{n} \in \Omega ,n=1,2,...)$ и последовательность целых чисел $(m_{n} )$ такие, что для любого действительного $t$ }

\[\mathop{\lim }\limits_{n\to \infty } F_{n} (s+it,\bar{\theta }_{n} )=\zeta (s+it),\]
\textit{равномерно в круге $|s-3/4|\le r$; здесь}

\[F_{n} (s+it,\bar{\theta }_{n} )=\prod _{p\le m_{n} } \left(1-e^{2\pi i\theta _{p}^{n} } p^{-s} \right)^{-1} ;\bar{\theta }_{n} =(\theta _{p}^{n} ),\]
\textit{компоненты $\bar{\theta }_{n} $ индексированы простыми числами  и  произведение взято по всем простым числам, удовлетворяющим требуемому неравенству.}
\end{theorem}

Следует отметить, что длина приближающего произведения зависит от $t$.
\begin{corollary}
\textit{. Гипотеза Римана верна, т. е. }
\[\zeta (s)\ne 0,\]
когда $\sigma >1/2$ .
\end{corollary}

\begin{center}

\textbf{2.  Дополнительные утверждения.}
\end{center}

\begin{lemm}
\textit{. Пусть ряд аналитических функций }

\[\sum _{n=1}^{\infty } f_{n} (s)\]

\textit{задан в односвязной области $G$ комплексной $s$-плоскости и абсолютно сходится почти всюду в $G$ в смысле Лебега и функция }

\[\Phi (\sigma ,t)=\sum _{n=1}^{\infty } |f_{n} (s)|\]
\textit{является суммируемой функцией в $G$. Тогда, данный ряд равномерно сходится в любой компактной подобласти области $G$; в частности сумма этого ряда аналитична в $G$.}
\end{lemm}

\textit{Доказательство}. Достаточно показать, что теорема верна для любой прямоугольной области в $G$. Пусть $C$ прямоугольник в $G$ и $C'$ другой прямоугольник внутри $C$, при этом, их стороны параллельны координатным осям. Мы можем предполагать, что на контуре этих прямоугольников ряд сходится почти всюду, в соответствии с теоремой Фубини (см. [7, стр.208]). Пусть $\Phi _{0} (s)=\Phi _{0} (\sigma ,t)$ является суммой данного ряда в точках сходимости По теореме Лебега об ограниченной сходимости (см. [21, стр.293]):

\[(2\pi i)^{-1} \int _{C} \frac{\Phi _{0} (s)}{s-\xi } ds=\sum _{n=1}^{\infty } (2\pi i)^{-1} \int _{C} \frac{f_{n} (s)}{s-\xi } ds,\]
где интегралы взяты в смысле Лебега. Поскольку на правой части равенства интегралы существуют также, в смысле Римана, то применяя формулу Коши, получаем

\[\Phi _{1} (\xi )=(2\pi i)^{-1} \int _{C} \frac{\Phi _{0} (s)}{s-\xi } ds=\sum _{n=1}^{\infty } f_{n} (s),\]
где $\Phi _{1} (\xi )=\Phi _{0} (\xi )$ почти всюду и $\xi $ произвольная точка на или внутри контура. Далее, ряд в $C'$ оценивается следующим образом

\[|f_{n} (\xi )|\le (2\pi )^{-1} \int _{C} \frac{|f_{n} (s)|}{|s-\xi |} |ds|\le (2\pi \delta )^{-1} \int _{C} |f_{n} (s)||ds|,\]
где $\delta $ обозначает минимальное расстояние между сторонами $C$ и $C'$. Ряд

\[\sum _{n=1}^{\infty } \int _{C} |f_{n} (s)||ds|\]
сходится в согласии с теоремой Лебега о монотонной сходимости (см. [21, стр.290]). Следовательно, ряд $\sum _{n=1}^{\infty } f_{n} (\xi )$ сходится равномерно внутри $C'$. Лемма 2 доказана.

Введем теперь понятие пространства Харди.
\begin{definition}
\textit{. Пространством Харди }\textbf{\textit{$H_{2}^{(R)} ,R>0$ }}\textit{называется множество функций $f(s)$, определенных для $\left|s\right|<R$ и аналитических в этой области, для которых}

\[\left\| f\right\| ^{2} =\mathop{\lim }\limits_{r\to R} \mathop{\int\!\!\!\!\int}\nolimits _{\left|s\right|<r}\left|f(s)\right| ^{2} d\sigma dt<\infty ;s=\sigma +it.   (*)\]
\end{definition}

(В настоящей статье мы придерживаемся терминологии кн. [24], называя пространство, введенное выше, пространством Харди. Заметим, что в книге [48] на стр. 177 введен класс \textbf{\textit{$H_{2}^{} $}} и имеется ссылка на книгу [49] где отмечено, что класс функций с интегральными средними

\[\mathop{\lim }\limits_{r\to 1} I_{r} =\mathop{\lim }\limits_{r\to 1} \int _{|z|=1}\left|f(z)\right|^{2} \left|dz\right|= \mathop{\lim }\limits_{r\to 1} r\int _{0}^{2\pi }\left|f(re^{i\theta } )\right|^{2} d\theta  <\infty \]
и его обобщения были исследованы Харди и Риссом. Этот класс оказывается замкнутой линейной оболочкой ортогональной на окружности \textit{$\left|z\right|=1$} системы \textit{$1,z,z^{2} ,...$}. На стр. 186 книги Уолша введен класс \textbf{\textit{$H'_{2} $}} , являюшийся замкнутой линейной оболочкой той же системы в области \textit{$C':\left|z\right|<1$} и оказывается, что это не что иное как лебеговый класс \textbf{\textit{$L_{2} $}} на \textit{$C':\left|z\right|<1$}, с точностью до эквивалентности в смысле Лебега. Здесь возникают средние (*). На стр. 189 той же книги отмечается, что Гика (1936 г.) и Бергман (1950 г.) исследовали подобные системы в многосвязных областях при помощи обобщенных средних вида (*)).

Очевидно, пространство Харди -- линейное пространство, в котором можно ввести скалярное произведение функций с помощью равенства
\begin{equation}\label{3}
\left(f(s),g(s)\right)=Re\mathop{\int\!\!\!\!\int}\nolimits_{\left|s\right|\le R}f(s)\overline{g(s)} d\sigma dt.
\end{equation}

Рассматривая \textbf{\textit{$H_{2}^{(R)} $}} как линейное пространство над полем вещественных чисел и используя введенное скалярное произведение, мы превращаем это пространство в вещественное гильбертово пространство.
\begin{lemm}
\textit{. Пространство }\textbf{\textit{$H_{2}^{(R)} $ }}\textit{вместе с введенным скалярным произведением (3) является вещественным гильбертовым пространством.}
\end{lemm}

\textit{ Доказательство}. Достаточно доказать, что произвольная фундаментальная последо-вательность $\left(f_{n} (s)\right)_{n\ge 1} $ сходится к некоторой аналитической функции . Поскольку, последовательность фундаментальная, то найдется такая последовательность натуральных чисел $\left(n_{k} \right)_{k\ge 0} $, что для любого натурального $k$

\[\left\| f_{n_{k} } -f_{n_{k-1} } \right\| \le 2^{-k} .\]

Рассмотрим ряд аналитических функций

\[f_{n_{0} } +\sum _{k=1}^{\infty }\left(f_{n_{k} } -f_{n_{k-1} } \right) \]
и докажем, что он сходится равномерно в любом замкнутом круге лежащем внутри круга $s<R$. Согласно (3) имеем:

\[\left\| f(s)\right\| ^{2} =\mathop{\int\!\!\!\!\int}\nolimits_{\left|s\right|<R}\left|f(s)\right|^{2} d\sigma dt. \]
Тогда обозначая $g(s)=\sum _{k=1}^{\infty }\left|f_{n_{k} } (s)-f_{n_{k-1} } (s)\right| $ получаем:

\[\mathop{\int\!\!\!\!\int}\nolimits_{\left|s\right|<R}g(s)d\sigma dt\le \sum _{k=1}^{\infty }\left(\pi R^{2} \mathop{\int\!\!\!\!\int}\nolimits_{\left|s\right|<R}\left|f_{n_{k} } -f_{n_{k-1} } \right|^{2} d\sigma dt \right)^{1/2}   \le \sqrt{\pi } R\sum _{k=1}^{\infty }2^{-k} <+\infty . \]
Следовательно, функция $g(s)$ суммируемая функция переменных $\sigma ,t$, и применима лемма 1. Применяя лемму 1, поучаем что ряд $f_{n_{0} } +\sum _{k=1}^{\infty }\left(f_{n_{k} } -f_{n_{k-1} } \right) $ сходится равномерно в любом круге  $\left|s\right|\le r<R$. Тогда подпоследовательность $\left(f_{n_{k} } (s)\right)_{k\ge 1} $ сходится к некоторой аналитической функции $\varphi (s)$. Так как последовательность фундаментальна, то для любого $\varepsilon >0$ найдется $n_{0} $ такое, что при любом натуральном $m>n_{0} $

\[\mathop{\int\!\!\!\!\int}\nolimits_{\left|s\right|<R}\left|\varphi (s)-f_{m} (s)\right|^{2} d\sigma dt<\varepsilon . \]

Пусть $r<R$ произвольное действительное число. Тогда используя соотношение кн.[23, стр. 345] получаем

\[r^{2} \left|\varphi (s)-f_{m} (s)\right|^{2} \le \pi ^{-1} \mathop{\int\!\!\!\!\int}\nolimits_{\left|s\right|<R}\left|\varphi (s)-f_{m} (s)\right|^{2} d\sigma dt<\varepsilon /\pi , \]

для любого $s,\left|s\right|\le r$. Поскольку $\varepsilon $ произвольно, то отсюда следует сходимость $\left(f_{n_{k} } (s)\right)_{k\ge 1} $. Следовательно рассматриваемое прстранство полное. Лемма 2 доказана.

Следующая лемма принадлежит С. М. Воронину ([28]) (мы приводим ее в слегка видоизмененной форме).
\begin{lemm}
\textit{. Пусть $g(s)$ является аналитической функцией в круге $|s|<r<1/4$, которая непрерывна и не обращается в нуль в замкнутом круге $|s|\le r$. Тогда для любого $\varepsilon >0$ и $y>2$ можно найти конечное множество простых чисел $M$, содержащее все простые $p\le y$ и некоторый элемент $\overline{\theta }=(\theta _{p} )_{p\in M} $ такой, что:}

1) $0\le \theta _{p} \le 1$ \textit{для} $p\in M$;

2) $\theta _{p} =\theta _{p}^{0} $ \textit{являются наперед заданными при} $p\le y$;

3) $\mathop{\max }\nolimits_{|s|\le r} |g(s)-\zeta _{M} (s+3/4;\overline{\theta })|\le \varepsilon ;$ \textit{здесь} $\zeta _{M} (s+3/4;\overline{\theta })$ \textit{определен равенством}

\[\zeta _{M} (s+3/4;\overline{\theta })=\prod _{p\in M} \left(1-e^{2\pi i\theta _{p} } p^{-s-3/4} \right)^{-1} .\]
\end{lemm}

\textit{ Доказательство. }Доказательство леммы 1 мы проведем следуя работе С. М. Воронина [28]. Поскольку \textit{$g(s)$}аналитична в круге мы рассмотрим вспомогательную функцию \textit{$g(s/\gamma ^{2} )$ }$(\gamma >1,\gamma ^{2} r<1/4)$, которая по непрерывности для любого $\varepsilon >0$ при некотором $\gamma $ удовлетвор-яет неравенству $\mathop{\max }\limits_{\left|s\right|\le r} \left|g(s)-g(s/\gamma ^{2} )\right|<\varepsilon $. Поэтому, достаточно доказать утверждение леммы для функции \textit{$g(s/\gamma ^{2} )$} в круге \textit{$|s|\le r$}. Преимущество состоит в том, что функция {$g(s/\gamma ^{2} )$} уже принадлежит пространству $H_{2}^{(\gamma r)} $ (круг имеет радиус, больше $r$, что важно для последующих рассуждений). Не нарушая, поэтому общности мы полагаем, что функция $g(s)$аналитична в круге \textit{$|s|\le r\gamma ^{2} $} и будем рассмотреть пространство $H_{2}^{(\gamma r)} $.

 Функция $\log g(s)$ по условиям теоремы не имеет внутри круга \textit{$|s|\le r\gamma $ }особенностей. Поэтому, достаточно доказать существование такого элемента \textit{$\overline{\theta }$}, удовлетворяющего условиям леммы 3, что

\[\mathop{\max }\nolimits_{|s|\le r} |g(s)-\zeta _{M} (s+3/4;\overline{\theta })|\le \varepsilon ;\]

Положим

\[u_{k} (s)=\log (1-e^{-2\pi i\theta _{k} } p_{k}^{-s-3/4} ),\]
беря главное значение логарифма. Используя разложение логарифмической функции в степенной ряд, мы можем писать

\[u_{k} (s)=-e^{-2\pi i\theta _{k} } p_{k}^{-s-3/4} +\nu (s),\quad \]

где
\[\left|\nu (s)\right|\le \left|(1/2)e^{-4\pi i\theta _{k} } p_{k}^{-2s-3/2} +\cdots \right|=O(p_{k}^{2r-3/2} ).\]
Поскольку $r<1/4$, мы можем найти $\delta >0$ такое, что $2\delta +2r-3/2<-1$. Тогда определение функции $u_{k} (s)$, вместе с последним неравенством показывает, что ряд
\begin{equation}\label{4}
\sum _{k=n+1}^{\infty } \eta _{k} (s);\eta _{k} (s)=-e^{-2\pi i\theta _{k} } p_{k}^{-s-3/4} ;n=\pi (y)
\end{equation}
отличается от ряда $\sum  u_{k} (s)$ на абсолютно сходящийся ряд. Сначала мы докажем, что для любой функции $\varphi (s)\in H_{2}^{(\gamma r)} $ пространства Харди , при подходящем $\bar{\theta }$, найдется некоторая перестановка ряда $\sum  \eta _{k} (s)$, сходящаяся по норме к функции $\varphi (s)$. Из этого будет следовать равномерная сходимость этого ряда, в круге \textit{$|s|\le r$,} к функции $\varphi (s)$, согласно лемме 1. В частности,  полагая

\[\varphi (s)=\log g(s)-\sum _{k>n}(u_{k} (s)-\eta _{k} (s))-\sum _{k\le n}u_{k} (s)  \]
и учитывая последнее замечание, мы можем найти некоторую престановку ряда $\sum _{k>n} \eta _{k} (s) $, сходящуюся к $\varphi (s)$. При этом, соотверствующая перестановка ряда $\sum _{k>n}(u_{k} (s)-\eta _{k} (s)) $ сходится к прежней своей сумме равномерно. Тогда для любого $\varepsilon $ найдется такое множество индексов $k\in M$, что

\[\mathop{\max }\limits_{\left|s\right|\le r} \left|\varphi (s)-\sum _{k\in M,\log p_{k} >y}\eta _{k} (s) \right|\le \varepsilon /2.\]

Пусть $q(s)=\sum _{k=n+1}^{\infty }\left(u_{k} (s)-\eta _{k} (s)\right) $. Поскольку, этот ряд сходится абсолютно, то  можно подоб-рать упомянутое множество $M$ так, чтобы выполнялось соотношение

\[\left|q(s)-\sum _{k\in M,k>n}\left(u_{k} (s)-\eta _{k} (s)\right) \right|\le \varepsilon /2.\]
Тогда получим:

\[\left|\varphi (s)-\sum _{k\in M,\log p_{k} >y}\eta _{k} (s) \right|=\left|\log g(s)-\sum _{n\in M}u_{n} (s) \right|\le \varepsilon ,\]
и, тем самым, доказательство леммы 3 будет завершено.

Рассмотрим ряд (\ref{4}) и применим теорему П. 6. 1.. Для этого докажем выполнимость условий этой теоремы при подходящим образом выбранном $\bar{\theta }$.

Сначала полагая $R=\gamma r$, рассмотрим пространство \textbf{\textit{$H_{2}^{(R)} $}}. Имеем:   $\left\| \eta _{k} (s)\right\| ^{2} =\mathop{\int\!\!\!\!\int}\nolimits_{\left|s\right|\le R}\left|e^{-2\pi i\theta _{k} } p_{k}^{-s-3/4} \right|^{2} d\sigma dt\le  \pi R^{2} p_{k}^{2r-3/2} .$
Следовательно,

\[\sum _{k=1}^{\infty }\left\| \eta _{k} (s)\right\| ^{2} \le \pi R^{2}  \sum _{k=1}^{\infty }p_{k}^{2r-3/2} <+\infty  ,\]
т. е. первое условие теоремы П. 6. 1. выполнено.

Пусть теперь$\varphi (s)\in H_{2}^{(R)} $произвольный элемент пространства с условием $\left\| \varphi (s)\right\| ^{2} =1.$ Пусть $\varphi (s)$ имеет следующий степенной ряд в круге $\left|s\right|\le R$:

\[\varphi (s)=\sum _{n=0}^{\infty }\alpha _{n} s^{n} . \]
Тогда,

\[1=\mathop{\int\!\!\!\!\int}\nolimits_{\left|s\right|\le R}\left|\sum _{n=0}^{\infty }\alpha _{n} s^{n}  \right| ^{2} d\sigma dt.\]
Произведем замену переменных под интегралом по формулам: $\sigma =r\cos \varphi ,\; t=r\sin \varphi ,\; r\le R$ $0\le \varphi <2\pi $. Тогда,

\[1=\sum _{n=0}^{\infty }\sum _{n=0}^{\infty }\alpha _{n} \bar{\alpha }_{m}   \int _{0}^{R}r^{n+m+1} \int _{0}^{2\pi }(\cos 2\pi (n-m)\varphi +i\sin 2\pi (n-m)\varphi )  d\varphi .\]
Внутренний интеграл равен 0, при $m\ne n$ и $2\pi $ в противном случае. Следовательно,
\begin{equation}\label{5}
\pi \sum _{n=0}^{\infty }\left|\alpha _{n} \right|^{2} R^{2n+2}  (n+1)^{-1} =1.]
\end{equation}

Докажем теперь, что существует не зависящая от функции $\varphi (s)$ точка $\bar{\theta }$, такая что ряд $\sum _{k=1}^{\infty }\left(\eta _{k} (s),\varphi (s)\right) $ сходится после некоторой перестановки его членов. Мы имеем

\[(\eta _{k} (s),\varphi (s))=-Re\int  \int _{|s|\le R} e^{-2\pi i\theta _{k} } p_{k}^{-s-3/4} \overline{\varphi (s)}d\sigma dt=Re[-e^{-2\pi i\theta _{k} } \Delta (\log p_{k} )],\]
где

\[\Delta (x)=\int  \int _{|s|\le R} e^{-x(s+3/4)} \overline{\varphi (s)}d\sigma dt.\]
$\Delta (x)$ можно представить в следующем виде:

\[\Delta (x)=e^{-3x/4} \mathop{\int\!\!\!\!\int}\nolimits_{\left|s\right|\le R}\left(\sum _{n=0}^{\infty }(-sx)^{n} /n! \right) \overline{\left(\sum _{n=0}^{\infty }\alpha _{n} s^{n}  \right)}d\sigma dt=\]

\[=\pi R^{2} e^{-3x/4} \sum _{n=0}^{\infty }(-1)^{n} \bar{\alpha }_{n} x^{n} R^{2n} /(n+1)! =\pi R^{2} e^{-3x/4} \sum _{n=0}^{\infty }\beta _{n} (xR)^{n} /n! ,\]
где $\beta _{n} =(-1)^{n} R^{n} \bar{\alpha }_{n} /(n+1)$. Из (\ref{5}) заключаем:

\[\sum _{n=1}^{\infty }\left|\beta _{n} \right| ^{2} \le 1.\]
Следовательно, $\left|\beta _{n} \right|\le 1$ и, поэтому функция
\begin{equation}\label{6}
F(u)=\sum _{m=0}^{\infty } \frac{\beta _{m} }{m!} u^{m}
\end{equation}
будет целой функцией, при этом

\[\Delta (x)=\pi R^{2} e^{-3x/4} F(xR).\quad \]

  Докажем, что для любого $\delta >0$ найдется стремящаяся к бесконечности последовательность $u_{1} ,u_{2} ,...$, удовлетворяющая неравенству
\begin{equation}\label{7}
|F(u_{j} )|>ce^{-(1+2\delta )u_{j} } .
\end{equation}
Допустим противное, т. е. пусть найдется положительное число $\delta <1$ такое, что при некотором, достаточно большом $A>0$ следующее неравенство

\[\left|F(u)\right|\le Ae^{-(1+2\delta )u} \]
выполнено для всех $u\ge 0$; в этом случае имеем:

\[\left|e^{(1+\delta )u} F(u)\right|\le Ae^{-\delta \left|u\right|} \; ;u\ge 0,\]
Из доказанного выше, при $u<0$ получаем:

\[\left|F(u)\right|\le \sum _{n=0}^{\infty }\left|u\right|^{n} /n!=e^{-u}  ,\]
из которого следует

\[\left|e^{(1+\delta )u} F(u)\right|\le e^{\delta \, u} \le e^{-\delta \left|u\right|} .\]
Из полученных оценок следует существование интеграла:

\[\int _{-\infty }^{\infty }\left|e^{(1+\delta )u} F(u)\right|^{2} du. \]
Поскольку, функция \ref{6} является целой функцией экспоненциального типа, то функция $e^{(1+\delta )u} F(u)$, также, будет таковым и принадлежит классу $E^{\sigma } $ (см [4, стр.408]), с $\sigma <3$.Тогда по теореме Пели -- Винера (см. там же) найдет финитная функция $f(\xi )\in L_{2} (-3,3)$ такая, что

\[e^{(1+\delta )u} F(u)=\int _{-3}^{3}f(\xi )e^{iu\xi } d\xi . \]
Взяв обратное преобразование Фурье находим:

\[f(\xi )=\frac{1}{2\pi } \int _{-\infty }^{\infty }\left(e^{(1+\delta )u} F(u)\right)e^{-iu\xi } du. \]

Из найденных выше оценок следует, что этот интеграл сходится абсолютно и равномерно в полосе $\left|Im\xi \right|<\delta /2$, и потому представляет собой аналитическую в этой полосе функции, что противоречит финитности $f(\xi )$. Полученное противоречие доказывает существование последовательности точек с условием (\ref{7}).

Обоначая $x_{j} =u_{j} /R$ на основании (\ref{7}) можем утверждать, что

\[\left|\Delta (x_{j} )\right|>ce^{-3x_{j} /4} \left|F(x_{j} R)\right|\ge ce^{-3x_{j} /4} e^{-(1+2\delta )x_{j} R} =ce^{-x_{j} (R+2\delta R+3/4)} .\]
Если $\delta >0$ достаточно мало, то $R+2\delta R+3/4<1$ и, следовательно, существует $\delta _{0} >0$ такое, что
\begin{equation}\label{8}
\left|\Delta (x_{j} )\right|>e^{-(1-\delta _{0} )x_{j} } .
\end{equation}
Рассмотрим функцию $\Delta (x)$ на отрезке [$x_{j} -1,x_{j} +1$]. Следуя [28], положим $N=\left[x_{j} \right]+1$. Из оценки для коэффициентов $\beta _{n} $ следует:

\[\left|\sum _{n=N^{2} +1}^{\infty }\frac{\beta _{n} }{n!} (xR)^{n}  \right|\le \sum _{n=N^{2} +1}^{\infty }\frac{(xR)^{n} }{n!}  \le \frac{(xR)^{N^{2} } }{(N^{2} )!} \sum _{n=0}^{\infty }\frac{(xR)^{n} }{n!} \le \frac{(xR)^{N^{2} } }{(N^{2} )!} e^{N} , \]
Так как, при целых $n,m\ge 0$ $(n+m)!=n!(n+1)\cdots (n+m)\ge n!m!.$ При достаточно больших натуральных $m$ по формуле Стирлинга имеем:

\[m!=\Gamma (m+1)\ge e^{m\log m-m} =(m/e)^{m} .\]
Следовательно,

\[\left|\sum _{n=N^{2} +1}^{\infty }\frac{\beta _{n} }{n!} (xR)^{n}  \right|\le \frac{(xR)^{N^{2} } }{(N^{2} )!} e^{N} \le N^{N^{2} } \left(\frac{N^{2} }{e} \right)^{N^{2} } e^{N} <<e^{-2x_{j} } ,\]
при $x\in $[$x_{j} -1,x_{j} +1$]. Далее, $\sum _{n=0}^{N^{2} }\beta _{n} (xR)^{n} /n!<<e^{xR}  $. Анологично,

\[\left|\sum _{n=N^{2} +1}^{\infty }\frac{(-3x/4)^{n} }{n!}  \right|\le \frac{(3x/4)^{N^{2} } }{(N^{2} )!} \sum _{n=0}^{\infty }\frac{(3x/4)^{n} }{n!} \le \frac{(3x/4)^{N^{2} } }{(N^{2} )!} e^{N} <<e^{-2x_{j} }  \]
и  $\sum _{n=0}^{N^{2} }(-3x/4)^{n} /n!<<e^{3x/4} , $ при $x\in $[$x_{j} -1,x_{j} +1$]. Таким образом,

\[\Delta (x)=\pi R^{2} \sum _{n=0}^{N^{2} }\frac{\left(-3x/4\right)^{n} }{n!}  \sum _{n=0}^{N^{2} }\frac{\beta _{n} }{n!} \left(xR\right)^{n}  +O(e^{-x_{j} } )=\sum _{n=0}^{N^{4} }a_{n} x^{n} +O(e^{-x_{j} } ). \]
Согласно (\ref{8}), мы получаем неравенство

\[\mathop{\max }\limits_{|x-x_{j} |\le 1} |\Delta (x)|>e^{-(1-\delta _{0} )x_{j} } ,\]
для любого$j=1,\, 2,\, ...$. Пусть $a_{n} =b_{n} +ic_{n} ,b_{n} ,c_{n} \in R$. Тогда,

\[\Delta (x)=\sum _{n=0}^{N^{4} }b_{n} x^{n} +i \sum _{n=0}^{N^{4} }c_{n} x^{n} +O(e^{x_{j} } ), \]
поэтому, для каждого$j$ хотя бы одно из следующих неравеств выполнено:

\[\mathop{\max }\limits_{|x-x_{j} |\le 1} \left|\sum _{n=0}^{N^{4} }b_{n} x^{n}  \right|>0.1e^{-(1-\delta _{0} )x_{j} } ,\]
или,

\[\mathop{\max }\limits_{|x-x_{j} |\le 1} \left|\sum _{n=0}^{N^{4} }c_{n} x^{n}  \right|>0.1e^{-(1-\delta _{0} )x_{j} } .\]

Рассмотрим сначала первую возможность. Пусть$x_{0} $-точка, где достигается максимум модуля. Обозначим $\tau _{j} $ отрезок в интервале [$x_{j} -1,x_{j} +1$] содержащий точку $x_{0} $, где имеет место упомянутая первая возможность. Пусть для определенности $g(x_{0} )<0;\, g(x)=\sum _{n=0}^{N^{4} }b_{n} x^{n}  $. Если $\tau _{j} $$\ne $[$x_{j} -1,x_{j} +1$] (случай совпадения интервалов тривиален), то найдется точка $x_{1} \in \tau _{j} $, для которой
\[\left|g(x_{1} )\right|\le 0.1\left|g(x_{0} )\right|.\]
Теперь имеем:

\[\left|g(x_{0} )-g(x_{1} )\right|\ge 0.5\left|g(x_{0} )\right|.\]
По теореме Лагранжа найдется точка $y_{j} \in \tau _{j} $ такая, что

\[\left|g'(y_{j} )(x_{1} -x_{0} )\right|\ge 0.5\left|g(x_{0} )\right|.\]
Применяя теорему П. 2. 9. кн. [24], находим:

\[N^{8} \left|g(x_{0} )\right|\left|x_{1} -x_{0} \right|\ge \left|g'(y_{j} )(x_{1} -x_{0} )\right|\ge 0.5\left|g(x_{0} )\right|.\]

Итак, отрезок $\tau _{j} $ имеет длину не меньше, чем $0.5x_{j}^{-8} $. Для определенности положим $\tau _{j} =\left[\alpha ,\alpha +\beta \right]$. По теореме Ш. Валле-Пуссена -- Ж. Адамара в интервале $\tau _{j} $ содержится

\[\int _{e^{\alpha } }^{e^{\alpha +\beta } }\frac{dx}{\log x}  +O(e^{\alpha +\beta } e^{c\sqrt{\alpha } } )=\int _{\alpha }^{\alpha +\beta }\frac{e^{u} }{u}  du+O(e^{\alpha +\beta } e^{c\sqrt{\alpha } } )\ge \]
\[\ge \frac{e^{\alpha } }{\alpha } \left[\left(e^{\beta } -1\right)+O\left(\frac{e^{\beta } }{e^{c\sqrt{\alpha } } } \right)\right]>>\frac{\beta e^{\alpha } }{\alpha } \]
простых чисел. Отберем те $p_{k} $ из них, для которых $p_{k} >y,k\equiv 0(mod4)$ и положим $\theta _{k} =0$. В случае, когда $g(x_{0} )>0$ будем полагать $\theta _{k} =1/2$ при $p_{k} >y,k\equiv 2(mod4)$. Тогда,

\[\sum _{\log p_{k} \in \tau _{j} ,k\equiv 0(mod4)}(\eta _{k} (s),\varphi (s)) =\sum _{\log p_{k} \in \tau _{j} ,k\equiv 0(mod4)}Re[-e^{-2\pi i\theta _{k} } \Delta (\log p_{k} )]>>\]
\[>> e^{x_{j} } e^{-(1-\delta _{0} )x_{j} } x_{j}^{-8} >> e^{\delta _{0} x_{j} /2}  .\]

Далее, при второй возможности, т. е. когда выполнено неравенство

\[\mathop{\max }\limits_{|x-x_{j} |\le 1} \left|\sum _{n=0}^{N^{4} }c_{n} x^{n}  \right|>0.1e^{-(1-\delta _{0} )x_{j} } \]
мы отбираем простые $p_{k} $, с условием $k\equiv 1(mod4)$, полагая, при этом, $\theta _{k} =1/4$, если значение многочлена в точке $x_{0} $ отрицательно; в противном случае берем $k$, с условием $k\equiv 3(mod4)$,  полагая, при этом $\theta _{k} =3/4$.

Таким образом, найдется бесконечное множество индексов с условием

\[\sum _{\log p_{k} \in \tau _{j} ,k\equiv 0\vee 2(mod4)}(\eta _{k} (s),\varphi (s)) >>e^{\delta _{0} x_{j} /2} ,\]
и бесконечное множество других значений $j$, для которых

\[-\sum _{\log p_{k} \in \tau _{j} ,k\equiv 1\vee 3(mod4)}(\eta _{k} (s),\varphi (s))>> e^{\delta _{0} x_{j} /2} .\]

Далее, из доказанных выше оценок мы заключаем, что

\[|\Delta (x)|\le \pi R^{2} e^{-x/2} ;\]
так, что $|(\eta _{k} (s),\varphi (s))|\to 0$, когда $k\to \infty $. Следовательно, ряд

\[\sum _{n=1}^{\infty } (\eta _{k} (s),\varphi (s))\]
содержит подряды, не имеющие общих компонент, расходящиеся, соответственно к $+\infty $ и к $-\infty $. Тогда некоторая перестановка ряда

\[\sum _{n=1}^{\infty } (\eta _{k} (s),\varphi (s))\]
сходится условно. Поэтому, по теоеме П. 6. 1., существует перестановка ряда $\sum _{p_{n} >y} u_{n} (s)$, сходящаяся к $\varphi (s)-\sum _{p_{n} \le y} u_{n} (s)$ равномерно. Взяв достаточно длинную частичную сумму, мы получаем нужный результат. Лемма 2 доказана.

\begin{center}

\textbf{3.  Основной вспомогательный результат.}

\end{center}

Пусть $\omega \in \Omega $, $\Sigma (\omega )=\{ \sigma \omega |\sigma \in \Sigma \} $ и $\Sigma '(\omega )$ означает замкнутое множество всех предельных точек последовательности $\Sigma (\omega )$. Для действительного $t$ мы обозначаем $\{ t\Lambda \} =(\{ t\lambda _{n} \} )$, где $\Lambda =(\lambda _{n} )$. Ниже мы обозначаем $\mu $ произведение линейных мер Лебега $m$ заданных на отрезке $[0,1]$: $\mu =m\times m\times \cdot \cdot \cdot $.
\begin{lemm}
\textit{. Пусть $A\subset \Omega $ является конечно- симметричным подмножеством меры нуль и $\Lambda =(\lambda _{n} )$ является неограниченной монотонно возрастающей последовательностью положительных действительных чисел, любое конечное подсемейство элементов которой линейно независимо над полем рациональных чисел. Пусть $B\supset A$ произвольное открытое подмножество с $\mu (B)<\varepsilon $ и }

\[E_{0} =\{ 0\le t\le 1|\{ t\Lambda \} \in A\wedge \Sigma '\{ t\Lambda \} \subset B\} .\]

\textit{Тогда, имеем $m(E_{0} )\le 6c\varepsilon $, где $c$ абсолютная постоянная и $m$ обозначает меру Лебега.}
\end{lemm}

\textit{Доказательство}. Пусть $\varepsilon $ является произвольным малым положительным числом. Поскольку числа $\lambda _{n} $ линейно независимы, то для любой конечной перестановки $\sigma $ имеем $(\{ t_{1} \lambda _{n} \} )\ne (\{ t_{2} \lambda _{\sigma (n)} \} )$, когда $t_{1} \ne t_{2} $. Действительно, в противном случае мы бы получили равенство $\{ t_{1} \lambda _{s} \} =\{ t_{2} \lambda _{s} \} $, для достаточно большого натурального $s$, т. е. $(t_{1} -t_{2} )\lambda _{s} =k,k\in Z$. Далее, записывая то же самое равенство для некоторого другого целого $r>m$ мы имеем соотношение

\[k_{1} /\lambda _{r} -k/\lambda _{s} =\frac{k_{1} \lambda _{s} -k\lambda _{r} }{\lambda _{r} \lambda _{s} } =0,\]
которое противоречит линейной независимости чисел $\lambda _{n} $. Следовательно для любой пары различных чисел $t_{1} $ и $t_{2} $ имеем $(\{ t_{1} \lambda _{n} \} )\notin \{ (\{ t_{2} \lambda _{\sigma (n)} \} )|\sigma \in \Sigma \} $. По условию леммы найдется семейство открытых шаров $B_{1} ,B_{2} ,...$ (в топологии Тихонова) такое, что каждый шар не содержит никакого другого шара из этого семейств (шар, содержащегося в другом может отброшен), при этом

\[A\subset B\subset \bigcup _{j=1}^{\infty } B_{j} ,\sum  \mu (B_{j} )<1.5\varepsilon .\]

Теперь берем некоторую перестановку $\sigma \in \Sigma $, определенную равенствами $\sigma (\ref{1})=n_{1} ,$ $...,\sigma (k)=n_{k} $, где натуральные числа подобраны следующим образом. Сначала мы берем $N$ такое, что

\[\mu (B'_{N} )<2\varepsilon _{1} ,\]
где $B'_{N} $ является проекцией шара $B_{1} $ в подпространство первых $N$ координатных осей и $\mu (B_{1} )=\varepsilon _{1} $. Покроем $B'_{N} $ кубами с ребром $\delta $ и суммарной мерой не превосходящей $3\varepsilon _{1} $. Положим $k=N$ и определим числа $n_{1} ,...,n_{k} $, используя следующие неравенства
\begin{equation}\label{9}
\lambda _{n_{1} } >1,\lambda _{n_{2} }^{-1} <(1/4)\delta \lambda _{n_{1} }^{-1} ,\lambda _{n_{3} }^{-1} <(1/4)\delta \lambda _{n_{2} }^{-1} ,...,\lambda _{n_{k} }^{-1} <(1/4)\delta \lambda _{n_{k-1} }^{-1} ,\delta <1.
\end{equation}

Теперь возьмем произвольный куб с ребром $\delta $ и с центром в некоторой точке $(\alpha _{m} )_{1\le m\le k} $. Тогда точка $(\{ t\lambda _{n_{m} } \} )$ будет принадлежать этому кубу, если

\[|\{ t\lambda _{n_{m} } \} -\alpha _{m} |\le \frac{\delta }{2} .\]
Из определения дробной части при $m=1$ для некоторого целого $r$ имеем:
\begin{equation}\label{10}
\frac{r+\alpha _{1} -\delta /2}{\lambda _{n_{1} } } \le t\le \frac{r+\alpha _{1} +\delta /2}{\lambda _{n_{1} } } .
\end{equation}
Мера множества таких $t$ не превосходит величину $\delta \lambda _{n_{1} }^{-1} $. Число же таких интервалов, соответствующих разным значениям $r=[t\lambda _{n_{1} } ]\le \lambda _{n_{1} } $ не превосходит

\[[\lambda _{n_{1} } ]+2\le \lambda _{n_{1} } +2.\]
Суммарная мера соответствующих интервалов не превосходит

\[\le (\lambda _{n_{1} } +2)\delta \lambda _{n_{1} }^{-1} \le (1+2\lambda _{n_{1} }^{-1} )\delta .\]

Теперь рассмотрим один из интервалов (\ref{10}); беря $m=2$, будем иметь
\begin{equation}\label{11}
\frac{s+\alpha _{2} -\delta /2}{\lambda _{n_{2} } } \le t\le \frac{s+\alpha _{2} +\delta /2}{\lambda _{n_{2} } } ,
\end{equation}
с $s=[t\lambda _{n_{2} } ]\le \lambda _{n_{2} } $. Поскольку мы рассматриваем условия (\ref{10}) и (\ref{11}) одновременно, то мы должны оценить суммарную меру тех интервалов (\ref{11}), которые имеют непустые пересечения с интервалами вида (\ref{10}), используя условия (\ref{9}). Число интервалов вида (\ref{11}) с длинами $\lambda _{n_{2} }^{-1} $, имеющих с одним интервалом вида (\ref{10}) непустое пересечение, не превосходит величину

\[[\delta \lambda _{n_{1} }^{-1} \lambda _{n_{2} } ]+2\le \delta \lambda _{n_{1} }^{-1} \lambda _{n_{2} } +2.\]
Тогда мера множества значений $t$, для которых одновременно выполнены условия (\ref{10}) и (\ref{11}) не превосходит

\[(\lambda _{n_{1} } +2)(2+\delta \lambda _{n_{1} }^{-1} \lambda _{n_{2} } )\delta \lambda _{n_{2} }^{-1} .\]

Можно продолжить эти рассуждения рассматривая все условия вида

\[\frac{l+\alpha -\delta /2}{\lambda _{n_{m} } } \le t\le \frac{l+\alpha +\delta /2}{\lambda _{n_{m} } } ,m=1,...,k.\]
Тогда находим следующую оценку для меры $m(\delta )$ множества тех $t$, для которых точки $(\{ t\lambda _{n_{m} } \} )$ попадают в рассматриваемые кубы с ребром $\delta $:

\[m(\delta )\le (2+\lambda _{n_{1} } )(2+\delta \lambda _{n_{1} }^{-1} \lambda _{n_{2} } )\ldots (2+\delta \lambda _{n_{k-1} }^{-1} \lambda _{n_{k} } ).\]

Проводя несложные преобразования, находим, учитывая условия (\ref{6}):

\[\le (2+\lambda _{n_{1} } )(2+\delta \lambda _{n_{1} }^{-1} \lambda _{n_{2} } )\cdots (2+\delta \lambda _{n_{k-1} }^{-1} \lambda _{n_{k} } )\delta \lambda _{n_{k} }^{-1} \le \delta ^{k} \prod _{m=1}^{\infty }(1+2m^{-2} ) \]
Суммируя по всем таким кубам, для меры множества тех $t$, для которых $(\{ t\lambda _{n_{m} } \} )\in B_{1} $, получаем верхнюю оценку вида $\le 3c\varepsilon ,c>0$. Заметим, что последовательность $\Lambda =(\lambda _{n} )$, определенная выше зависит от $\delta $. Мы будем для каждого шара $B_{k} $ фиксировать некоторую последовательность $\Lambda _{k} $, используя условия (\ref{9}). Рассматривая все такие шары мы обозначаем $\Sigma _{0} =\{ \Lambda _{k} |k=1,2,...\} $ . Так как множество $A$ является конечно-симметричным, то мера интересующего нас множества значений $t$ можно оценить используя произвольную после-довательность $\Lambda _{k} $, т. к., как было показано выше, множества $\Sigma (\{ t\Lambda \} )$ для разных значений $t$ имеют пустое пересечение.

 Докажем, что для любой точки $t\in E_{0} $, множество $\Sigma (\{ t\Lambda \} )$ содержится в объединении $\bigcup _{k\le n}B_{k}  $, при некотором $n$. Действительно, пусть для некоторого $t\in E_{0} $ все члены последовательности $\Sigma (\{ t\Lambda \} )$ не содержатся в объединении $\bigcup _{k\le n}B_{k}  $, для любого натурального $n$. Возможны два случая: 1) найдется точка $\bar{\theta }\in \Sigma (\{ t\Lambda \} )$, принадлежащая бесконечному множеству шаров $B_{k} $; 2) Найдется последовательность элементов $\bar{\theta }_{j} ,\bar{\theta }_{j} \in \Sigma (\{ t\Lambda \} )$, которая не содержится в произвольном конечном объединении шаров $B_{k} $. Рассмотрим обе возможности в отдельности и докажем, что они приводят к противоречию.

1) Пусть $\bar{\theta }\in B_{k_{1} } $. $B_{k_{2} } ,B_{k_{3} } ,...$ все шары, к которым принадлежит элемент $\bar{\theta }$. Обозначим $d$ расстояние от $\bar{\theta }$ до границы $B_{k_{1} } $. Поскольку $B_{k_{1} } $- открытое множество, то $d>0$. Пусть  $B_{k} $ произвольный шар радиуса $<d/2$, из числа содержащих точку $\bar{\theta }$. Из сказанного следует, что шар $B_{k} $ должен содержатся внутри шара $B_{k_{1} } $. Но это противоречит принятому выше соглашению.

2) Пусть $\bar{\theta }$ некоторая предельная точка последовательности $\bar{\theta }_{j} $. Согласно условию $\bar{\theta }\in B_{s} $ при некотором $s$. Пусть $d$ расстояние от $\bar{\theta }$ до границы $B_{s} $. Поскольку $\bar{\theta }$ предельная точка, то шар с центром в точке $\bar{\theta }$ и срадиусом $d/4$ содержит бесконечное множество членов последовательности $\bar{\theta }_{j} $, скажем члены $\bar{\theta }_{j_{1} } ,\theta _{j_{2} } ,...$. Так как, по 1), каждая точка этой последовательноси может принадлежать только конечному числу шаров, то указанная последовательность содержится в объединении бесконечного подсемейства шаров $B_{k} $. Среди них имеются бесконечное количество шаров, имеющих радиус $<d/4$. Все они, тогда должны содержатся внутри шара $B_{s} $. Полученное противоречие исключает и случай 2).

Итак, для любого $t\in E_{0} $ найдется $n$ такое, $\Sigma (\{ t\Lambda \} )\subset $$\bigcup _{k\le n}B_{k}  $. Отсюда, в свою очередь следует, что множество $E_{0} $ может быть представлено в виде объединения подмножеств $E_{k} ,k=1,2,{\rm ...}$, где $E_{k} =\{ t\in E_{0} |\Sigma (t\Lambda )\subset \bigcup _{s<k} B_{s} \} .$Тогда,

\[E_{k} =\{ t\in E_{0} |\Sigma (t\Lambda )\subset \bigcup _{m\le k} B_{m} \} ,E_{0} =\bigcup _{k=1}^{\infty } E_{k} ;E_{k} \subset E_{k+1} (k\ge 1).\]

Далее, согласно [23, стр.368], $m(E_{0} )=\mathop{\lim }\limits_{k\to \infty } m(E_{k} )$. Как было отмечено,$m(E_{k} )$ можно оценить пользуясь произвольной последовательностью $\Lambda '$:

\[m(E_{k} )=\mathop{\lim \sup }\limits_{\Lambda '\in \Sigma _{0} } m(E_{k} (\Lambda ')),\]
где $E_{k} (\Lambda ')=\{ t\in E_{k} |(\{ t\Lambda '\} )\in \bigcup _{s\le k}B_{s}  \} $. Следовательно,

\[m(E_{k} (\Lambda '))\le \sum _{s\le k}m (E^{(s)} (\Lambda ')),\]
где $E^{(k)} (\Lambda ')=\{ t\in E_{0} |(\{ t\Lambda '\} )\in B_{k} \} $. Применяя лемму 3 находим (выбирая подходящим образом $\Lambda '$):

\[m(E(\Lambda '))\le 6c(\varepsilon _{1} +\cdots +\varepsilon _{k} ).\]
Переходя к пределу, при $k\to \infty $, получаем требуемое соотношение. Доказательство леммы 3 завершено.

\begin{center}

\textbf{4.  Локальное приближение.}

\end{center}

\begin{lemm}
\textit{. Существует последовательность точек $(\bar{\theta }_{k} )$ $(\bar{\theta }_{k} \in \Omega )$ и натуральных чисел }

\textit{$(m_{k} )$ такие, что $\bar{\theta }_{k} \to 0$ и }

\[\mathop{\lim }\limits_{k\to \infty } F_{k} (s+3/4,\bar{\theta }_{k} )=\zeta (s+3/4)\]

\textit{в круге $|s|\le r,0<r<1/4$ равномерно по $s$.}
\end{lemm}

\textit{Доказательство.} Пусть $y>2$ является положительным целым числом, точнее определяемым ниже. Мы полагаем

\[y_{0} =y,y_{1} =2y_{0} ,...,y_{m} =2y_{m-1} =2^{m} y_{0} ,....\]
Из леммы 1 следует, что для данного положительного $\varepsilon $ и числа $y>2$ найдется множество простых $M_{1} $ и точка $\bar{\theta }_{1} =(\theta _{p}^{0} )_{p\in M_{1} } $ такая, что $M_{1} $ содержит все простые $p\le y$, с $\theta _{p}^{0} =0$ и

\[\mathop{\max }\limits_{|s|\le r} |\zeta (s+3/4)-\eta _{1} (s+3/4)|\le \varepsilon ;\eta _{1} (s+3/4)=\prod _{p\in M_{1} } \left(1-e^{2\pi i\theta _{p}^{0} } p^{-s-3/4} \right)^{-1} .\]
Теперь, обозначая $m_{1} =\mathop{\max }\nolimits_{m\in M_{1} } m$, положим

\[F_{1} (s+3/4;\bar{\theta })=\prod _{p\le m_{1} } \left(1-e^{2\pi i\theta _{p}^{0} } p^{-s-3/4} \right)^{-1} \]
и

\[h_{1} (s+3/4;\bar{\theta })=F_{1} (s+3/4;\bar{\theta })\prod _{p\in M_{1} } \left(1-e^{2\pi i\theta _{p}^{0} } p^{-s-3/4} \right)-1;\]
здесь $\theta _{p} =\theta _{p}^{0} $ при $p\in M_{1} $. Пусть $n$ обозначает натуральное число, каноническое разложение которого содержит только простые $p$, $p\bar{\in }M_{1} ,p\le m_{1} $, и

\[a_{n} (\bar{\theta })=e^{2\pi i\sum _{p\backslash n}\alpha _{p} \theta _{p}  } ;n=\prod  p^{\alpha _{p} } ;\]
Если $r+\delta <1/4$, то мы имеем

\[\int _{\Omega _{1} } \left(\int  \int _{|s|\le r+\delta } |h_{1} (s+3/4;\bar{\theta })|^{2} d\sigma dt\right)d\bar{\theta }\le \]

\[\le \int  \int _{|s|\le r+\delta } \left(\int _{\Omega _{1} } |h_{1} (s+3/4;\bar{\theta })|^{2} d\bar{\theta }\right)d\sigma dt\le \]

\[\le \pi (r+\delta )^{2} \mathop{\max }\limits_{|s|\le r+\delta } \int _{\Omega _{1} } |\sum _{n>y} a_{n} (\bar{\theta })n^{-s+3/4} |^{2} d\bar{\theta }\le \frac{4\pi (r+\delta )^{2} }{1-4r-4\delta } y^{-1/2+2r+2\delta } ;\]
здесь $\Omega _{1} $ означает проекцию $\Omega $ в подпространство координатных осей $\theta _{p} ,p\bar{\in }M_{1} $. Тогда из неравенства, полученного выше следует существование $\bar{\theta '}_{1} =(\theta _{p} )_{p\overline{\in }M_{1} } $ такого, что

\[\int  \int _{|s|\le r+\delta } |h_{1} (s+3/4;\bar{\theta '}_{1} )|^{2} d\sigma dt\le \frac{4\pi (r+\delta )^{2} }{1-4r-4\delta } y^{-1/2+2r+2\delta } \]
или

\[\mathop{\max }\limits_{|s|\le r} |h_{1} (s+3/4;\bar{\theta '}_{1} )|\le \sqrt{2} \delta ^{-1} \left(\frac{1}{2\pi } \int  \int _{|s|\le r} |h_{1} (s+3/4;\bar{\theta '}_{1} )|^{2} d\sigma dt\right)^{1/2} \le\]
\[\le c(\delta )y^{\delta +r-1/4} \]
(см. [22, стр. 345]) с постоянной $c(\delta )>0$. Следовательно, взяв $\bar{\theta }_{1} =(\bar{\theta }_{0} ,\bar{\theta '}_{1} ),$ где $\bar{\theta }_{0} =(\theta _{p}^{0} )_{p\in M_{1} } $, а $y=y_{0} $- удовлетворяющим условию

\[(A+1)c(\delta )y_{0}^{r+\delta -1/4} \le \varepsilon ;A=\mathop{\max }\limits_{|s|\le r} |\zeta (3/4+s)|,\]
будем иметь

\[\mathop{\max }\limits_{|s|\le r} \left\{|\zeta (3/4+s)-F_{1} (3/4+s;\bar{\theta }_{1} )|\right\}\le \]

\[\le \mathop{\max }\limits_{|s|\le r} \left\{|\zeta (3/4+s)-\eta _{1} (3/4+s)|+|\eta _{1} (3/4+s)|\cdot |h_{1} (3/4+s;\bar{\theta '}_{1} )|\right\}\le \]

\[\le \varepsilon +(A+1)c(\delta )y_{0}^{r+\delta -1/4} \le 2\varepsilon ,\]
Теперь, заменим $\varepsilon $ на $\varepsilon /2$. Существует множество $M_{2} $ простых чисел, содержащее все простые числа $p\le 2y_{0} =y_{1} $ и удовлетворяющее, по лемме 1, условиям

\[\mathop{\max }\limits_{|s|\le r} |\zeta (3/4+s)-\eta _{2} (3/4+s)|\le \varepsilon /2,\]
где

\[\eta _{2} (s+3/4)=\prod _{p\in M_{2} } \left(1-e^{2\pi i\theta _{p}^{(1)} } p^{-s-3/4} \right)^{-1} \]
и $\theta _{p}^{(\ref{1})} =0$ при $p\le y_{1} $. Анологичным образом, как выше, определяем функции

\[F_{2} (s+3/4;\bar{\theta })=\prod _{p\le m_{2} } \left(1-e^{2\pi i\theta _{p} } p^{-s-3/4} \right)^{-1} ;m_{2} =\mathop{\max }\limits_{m\in M_{2} } m\]
и
\[h_{2} (s+3/4;\bar{\theta })=F_{2} (s+3/4;\bar{\theta })\prod _{p\in M_{1} } \left(1-e^{2\pi i\theta _{p} } p^{-s-3/4} \right)-1;\]

Подобным образом, мы находим $\bar{\theta '}_{2} \in \Omega _{2} $ ($\Omega _{2} $ является проекцией $\Omega $ в подпространство координат $\theta _{p} ,p\overline{\in }M_{2} $) такой, что

\[\mathop{\max }\limits_{|s|\le r} |\zeta (3/4+s)-F_{2} (3/4+s;\bar{\theta }_{2} )|\le 2^{1+(r+\delta -1/4)} \varepsilon ,\bar{\theta }_{2} =(\bar{\theta }_{1} ,\bar{\theta '}_{2} ).\]
Действительно,

\[|F_{2} (3/4+s)-\eta _{2} (3/4+s)|=|\eta _{2} (3/4+s)|\cdot |h_{2} (3/4+s;\bar{\theta '}_{2} )|.\]
Теперь, взяв средние значения, как выше, получаем

\[\mathop{\max }\limits_{|s|\le r} |h_{2} (s+3/4;\bar{\theta '}_{2} )|\le \sqrt{2} \delta ^{-1} \left(\frac{1}{2\pi } \int  \int _{|s|\le r} |h_{2} (s+3/4;\bar{\theta '}_{2} )|^{2} d\sigma dt\right)^{1/2} \le \]
\[\le c(\delta )(2y_{0} )^{\delta +r-1/4} .\]
Следовательно,

\[\mathop{\max }\limits_{|s|\le r} |\zeta (3/4+s)-F_{2} (3/4+s;\theta _{2} )|\le \varepsilon /2+2^{1+(r+\delta -1/4)} \varepsilon ,\bar{\theta }_{2} =(\theta _{1} ,\bar{\theta '}_{2} ).\]

Повторяя рассуждения, для каждого $k>1$ находим $\bar{\theta }_{k+1} =(\bar{\theta }_{k} ,\bar{\theta '}_{k+1} )\in \Omega ,\bar{\theta }_{k} =(\bar{\theta }_{p}^{k} )_{p\in M_{k+1} } $ такой, что $\theta _{p}^{k} =0$, когда $p\le y_{k} $ и

\[\mathop{\max }\limits_{|s|\le r} |\zeta (3/4+s)-F_{k+1} (3/4+s;\bar{\theta }_{k+1} )|\le 2^{1+k(r+\delta -1/4)} \varepsilon ;\]
здесь

\[F_{k+1} (s+3/4;\bar{\theta })=\prod _{p\le m_{k+1} } \left(1-e^{2\pi i\theta _{p}^{0} } p^{-s-3/4} \right)^{-1} ;m_{k+1} =\mathop{\max }\limits_{m\in M_{k+1} } m.\]
Поэтому, равномерно по $s,|s|\le r$ имеем

\[\mathop{\lim }\limits_{k\to \infty } F_{k} (3/4+s,\bar{\theta }_{k} )=\zeta (3/4+s).\]
Лемма 4 доказана.

\begin{center}
\textbf{5.  Доказательство теоремы.}
\end{center}

Теперь рассмотрим интеграл

\[B_{k} =\int _{\Omega } \left(\int  \int _{|s|\le r} |F_{k+1} (3/4+s;\bar{\theta }_{k+1} +\bar{\theta })-F_{k} (3/4+s;\bar{\theta }_{k} +\bar{\theta })|d\sigma d\tau \right)d\bar{\theta },\]
где $k=0,1,...$, и если $k=0$, то полагаем $F_{0} (3/4+s;\bar{\theta }_{0} +\bar{\theta })=0$. Применяя неравенство Шварца и меняя порядок интегрирования, находим, как выше:

\[B_{k}^{2} \le 4\pi r^{2} \int  \int _{|s|\le r} d\sigma d\tau \int _{\Omega } |\prod _{p\le 2^{k-1} y_{0} } \left(1-e^{-2\pi i(\theta _{p}^{n} +\theta _{p} )} p^{-s-3/4} \right)^{-1} |^{2} \prod _{p\le 2^{k-1} y_{0} } d\theta _{p} \times \]

\[\times \sum _{n>2^{k-1} y_{0} } n^{2r+2\delta -3/2} \le c_{\delta } \left(2^{k-1} y_{0} \right)^{2r+2\delta +1-1/2} ;c_{\delta } >0.\]
Поскольку $2r+2\delta -1/2<0$, то из этой оценки следует сходимость почти всюду по $\bar{\theta }$ нижеприведенного ряда (для всех $\bar{\theta }\in \Omega _{0} $, где $\Omega _{0} $ является подмножеством полной меры и множество $A=\Omega \backslash \Omega _{0} $ конечно симметрично)
\begin{equation}\label{12}
\sum _{k=1}^{\infty } \int  \int _{|s|\le r} |F_{k} (3/4+s;\bar{\theta }_{k} +\bar{\theta })-F_{k-1} (3/4+s;\bar{\theta }_{k-1} +\bar{\theta })|d\sigma d\tau .
\end{equation}

По теореме Егорова (см. [10, стр.166]) этот ряд сходится почти равномерно вне некоторого множества $\Omega '_{1} ,\mu (\Omega '_{1} )=0$. Мы можем предполагать, что множество $A\bigcup \Omega '_{1} $ конечно симметрично (в противном случае можно брать множество все конечных перестановок всех его элементов). Найдется некоторое счетное семейство шаров $B_{r} $ с суммарной мерой не превосходящей $\varepsilon $, объединение которых содержит множество $A\bigcup \Omega '_{1} $. Для каждого нату-рального \textit{n} мы определяем множество $\Sigma '_{n} (t\Lambda )$ как множество всех предельных точек последовательности $\Sigma _{n} (\bar{\omega })=\{ \sigma \bar{\omega }|\sigma \in \Sigma \wedge \sigma (\ref{1})=1\wedge \cdots \wedge \sigma (n)=n\} $. Пусть $B^{(n)} =\{ t|\{ t\Lambda \} \in A\wedge $ $\sum '_{n} (\{ t\Lambda \} )\subset \bigcup _{r=1}^{\infty }B_{r}  \} ,\lambda _{n} =(1/2\pi )\log p_{n} $,$n=1,2,...$.При каждом $t$ последовательность $\sum _{n+1} (\{ t\Lambda \} )$ является подпоследовательностью последовательности $\sum _{n} (\{ t\Lambda \} )$. Поэтому, $\sum '_{n+1} (\{ t\Lambda \} )\subset $ $\subset \sum '_{n} (\{ t\Lambda \} )$. Следовательно, имеем $B^{(n)} \subset B^{(n+1)} $. Тогда, если обозначим $B=\bigcup _{n} B^{(n)} $, то получим $m(B)\le \sup m(B^{(n)} )$.

Оценим  $m(B^{(n)} )$. Множество $\sum '_{n} (\{ t\Lambda \} )$ является замкнутым множеством. Ясно, что, если мы будем "укоротить" последовательности $\{ t\Lambda \} $, оставляя только компоненты $\{ t\lambda _{n} \} $ с индексами, большими, чем $n$ и обозначим укороченную последовательность как $\{ t\Lambda \} '\in \Omega $, то множество $\sum '(\{ t\Lambda \} ')$ также будет замкнутым. Теперь рассмотрим произведения $[0,1]^{n} \times \{ \{ t\Lambda \} '\} $ (внешние скобки обозначают одноэлементное множество) для каждого $t$. Имеем

\[\{ t\Lambda \} \in [0,1]^{n} \times \{ \{ t\Lambda \} '\} \subset A.\]

Нижеследующий пример показывает, что из выполнимости последнего соотношения не следует равенство $A=\Omega $. Пусть $I=[0,1];U=[0;1/2];V=[1/2;1]$ и

\[X_{0} =U\times U\times \ldots ,X_{1} =V\times U\times \ldots ,\]

\[X_{2} =I\times V\times U\times \ldots ,X_{s+1} =I^{s} \times V\times U\times \ldots ,....\]
Ясно, что $\mu (X_{s} )=0$ для всех $s$. Тогда $\mu (X)=0$, где

\[X=\bigcup _{s=0}^{\infty } X_{s} .\]
Как видно из конструкции $X$, равенство

\[X=[0,1]^{s} \times X\]
справедливо для любого натурального $s$. Поскольку множество $[0,1]^{n} \times \{ \{ t\Lambda \} '\} $ замкнуто, то существует только конечное множество $R$ натуральных чисел такое, что $[0,1]^{n} \times \{ \{ t\Lambda \} '\} \subset \bigcup _{r\in R} B_{r} $. Рассмотрим множество "укороченных" точек $\bar{\theta }'$ шаров $B_{r} $. Пусть $B'_{r} =\{ \bar{\theta }'|\bar{\theta }\in B_{r} \} $. Из предыдущего соотношения следует $\{ t\Lambda \} '\in B'_{r} $ для всех $r\in R$. Тогда пересечение $\bigcap _{r\in R}B'_{r}  $, будучи открытым множеством, содержит точку $\{ t\Lambda \} '$. Внутри пересечения $\bigcap _{r\in R}B'_{r}  $ всегда найдется шар $K$ с центром в точке $\{ t\Lambda \} '$ и достаточно малого радиуса такой, что

\begin{equation}\label{13}
[0,1]^{n} \times \{ \{ t\Lambda \} '\} \subset [0,1]^{n} \times K\subset \bigcup _{r\in R} B_{r} ,
\end{equation}
для каждой рассматриваемой точки $t$. Аналогичное соотношение справедливо и в случае, когда точка $\{ t\Lambda \} $ будет заменена произвольной предельной точкой $\bar{\omega }$ последовательности $\Sigma (\{ t\Lambda \} )$, потому, что $\bar{\omega }\in B_{r} $. Если через $B'$ обозначить объединение всех открытых шаров  $K$, отвечающих всевозможным значениям $t$ и предельной точки $\bar{\omega }$, то получим:

\[\{ t\Lambda \} \in [0,1]^{n} \times \{ \{ t\Lambda \} '\} \subset A\subset [0,1]^{n} \times B'\subset \bigcup _{r=1}^{\infty }B_{r}  ,\]
для каждого рассматриваемого значения $t$, или

\[\{ \bar{\omega }\} \in [0,1]^{n} \times \{ \bar{\omega }\} '\subset A\subset [0,1]^{n} \times B'\subset \bigcup _{r=1}^{\infty }B_{r}  ,\]
для каждой предельной точки $\bar{\omega }$. Из этого следует, что $\mu ^{*} (B')\le \varepsilon $, где $\mu ^{*} $ означает внешнюю меру. Множество $B'$ является открытым и $\Sigma '(\{ t\Lambda \} ')\in B'$. Теперь мы можем применить лемму 3 и получить оценку $m(B^{(n)} )\le 6c\varepsilon $. Таким образом, имеем $m(B)\le 6c\varepsilon $.

 Пусть $t\notin B$. Тогда, $t\notin B^{(n)} $, для каждого $n=y_{k} ,k=1,2,3,...$. При каждом $k$, существует такая предельная точка $\bar{\omega }_{k} \in \Omega \backslash \bigcup _{r} B_{r} $ последовательности $\sum _{n} (\{ t\Lambda \} )$, для которой ряд

\[\sum _{l=1}^{\infty } \int  \int _{|s|\le r} |F_{l} (3/4+s;\bar{\theta }_{l} +\bar{\omega }_{k} )-F_{l-1} (3/4+s;\bar{\theta }_{l-1} +\bar{\omega }_{k} )|d\sigma d\tau \]
сходится. Так как множество $\Omega \backslash \bigcup _{r} B_{r} $ замкнуто, то предельная точка $\overline{\omega }=(\{ t\Lambda \} )$ последова-тельности $(\bar{\omega }_{k} )$ будет принадлежать множеству $\Omega \backslash \bigcup _{r} B_{r} $. Поэтому, ряд

\[\sum _{l=1}^{\infty } \int  \int _{|s|\le r} |F_{l} (3/4+s;\bar{\theta }_{l} +i\{ t\Lambda \} )-F_{l-1} (3/4+s;\bar{\theta }_{l-1} +i\{ t\Lambda \} )|d\sigma d\tau .\quad (14)\]
сходится, так как он сходится на этом множестве равномерно. Итак последний ряд сходится для всех $t$, за исключением значений $t$ из некоторого множества меры не превосходящей $12c\varepsilon $. В силу произвольности $\varepsilon $, последний результат показывает сходимость ряда (14) для почти всех $t$ (ясно, что условие $0\le t\le 1$ теперь может быть опущено). Тогда, по лемме 2, для произвольно взятого $\delta _{0} <1$ последовательность

\[F_{k} (3/4+s;\bar{\theta }_{k} +i\{ t\Lambda \} ),\quad (15)\]
для всех таких $t$ сходится в круге $|s|\le r\delta _{0} (\delta _{0} <1)$, равномерно, к некоторой аналитической функции $f(s+3/4;t)$:

\[\mathop{\lim }\limits_{k\to \infty } F_{k} (3/4+s+it;\bar{\theta }_{k} )=f(s+3/4;t).\]

Несмотря на полученный результат мы не можем использовать $t$ как переменную, так как левая и правая части этого равенства могут отличаться своими аргументами (правая часть определена как предел последовательности (15), где $t$ входит в выражение разрывной функции). Следовательно, принцип аналитического продолжения нельзя применить. Для того, чтобы завершить доказательство теоремы возьмем произвольное большое действительное число $T$. Поскольку, рассматриваемые значения $t$ всюду плотны в интервале$[-T,T]$, объединение кругов $C(t)=\{ 3/4+it+s:|s|\le r\delta _{0} \} $ содержит прямоугольник $3/4-r\delta _{0}^{2} \le Re(s+3/4)\le 3/4+r\delta _{0}^{2} ,-T\le Im(s+3/4)\le T$, в котором условия леммы 2 выполнены для ряда

\[F_{1} (s+3/4;\bar{\theta }_{1} )+(F_{2} (s+3/4;\bar{\theta }_{2} )-F_{1} (s+3/4;\bar{\theta }_{1} ))+\ldots .\quad (16)\]
Следовательно, по лемме 2, этот ряд определяет аналитическую функцию внутри рассматриваемого прямоугольника, которая совпадает с $\zeta (3/4+s)$ внутри круга $C(0)$. Для того, чтобы применить принцип аналитического продолжения возьмем односвязную открытую область, где обе функции $\log F_{*} (s)$ и $\log \zeta (s)$ регулярны (здесь функция $F_{*} (s)$ является суммой ряда (16)). Пусть $\rho _{1} ,...,\rho _{L} $ обозначают все возможные нули функции $\zeta (s)$ в рассматриваемом прямоугольнике, контур которого не содержит нулей функции $\zeta (s)$. Проведем разрезы вдоль отрезков $1/2\le Res\le Re\rho _{l} ,Ims=Im\rho _{l} ,l=1,...,L.$ В открытой области рассматриваемого прямоугольника, не содержащей указанные отрезки, функции $\log F_{*} (s)$ и $\log \zeta (s)$ регулярны. Тогда, равенство $F_{*} (s)=\zeta (s)$ выполнено во всей открытой области, определенной выше. Теперь мы получаем справедливость соотношения $F_{*} (s)=\zeta (s)$ во всем прямоугольнике, где обе функции регулярны. Доказательство теоремы завершено.

\begin{center}

\textit{}

\textbf{6.  Доказательство следствия.}

\end{center}

Вывод следствия получается применением теоремы Руше (см. [19, стр.137]). Пусть $t$ произвольное действительное число. Докажем, что для любого $0<r<3/4$ внутри круга $C=\{ s||s-3/4-it|=r\} $ $\zeta (s)$ не имеет нулей. Пусть радиус $r$ выбран так, что окружность $C$ не содержит нулей $\zeta (s)$, и мы полагаем

\[m=\mathop{\min }\limits_{s\in C} |\zeta (s)|.\]
По теореме найдется $n=n(t)$ такое, что на и внутри контура $C$ следующее неравенство выполнено

\[|\zeta (s)-F_{n} (s;\bar{\theta }_{n} )|\le 0.25m.\]
Тогда, на контуре $C$ верно неравенство:

\[|\zeta (s)-F_{n} (s;\bar{\theta }_{n} )|<|\zeta (s)|\]
По теореме Руше функции $\zeta (s)$ и $F_{n} (s;\bar{\theta }_{n} )$ имеют одинаковое число нулей внутри $C$. Но, функция $F_{n} (s;\bar{\theta }_{n} )$ не имеет нулей внутри круга $C$. Следовательно, $\zeta (s)$ также не имеет нулей внутри круга $C$. Так как, $t$ произвольно, то из последнего мы заключаем, что полоса $-r<Re\, s-3/4<r$ (для любого $0<r<1/4$) свободна от нулей функции $\zeta (s)$. Следствие доказано.

\begin{center}
\textbf{Литература}
\end{center}
\small
\begin{enumerate}
\setlength{\itemsep}{-1.5mm}

\item B. Bagchi. A joint universality theorem for Dirichlet L -functions. Math. Zeit., 1982,v. 181, p.319-335.

\item H. Bohr, R.Courant. Neue Anwendungen der Theorie der Diophantischen Approxima-tionen auf die Riemannsche Zetafunction.
 J.reine angew. Math. (1944) p.249-274.

\item K.Чандрасекхаран. Введение в аналитическую теорию чисел. M.: Мир, 1974.

\item А. Зигмунд. Тригонометрические ряды. т. 2., М.: Мир, 1965.

\item Р.Курант. Дифференциальное и интегральное исчисление. M. Наука, 1967,510 стр.

\item Г. Дэвенпорт. Мултипликативная теория чисел.-М. Наука, 1971.

\item Данфорд Н. и Шварц Дж.Е. Линейные операторы. Общая теория. М.ИИЛ.1962, 896 стр.

\item K.Чандрасекхаран. Арифметические функции. M. Наука, 1975., 270 стр.

\item L.Euler. Введение в анализ бесконечно малых. - M. :ОНТИ, 1936.

\item E. Хьюитт и K.A.Росс. Абстрактный гармонический анализ.т.I, Наука, 1975, 656 стр.

\item И. Ш. Джаббаров. О некоторых рядах Дирихле. MAK 98, Материалы первой региональной конференции, Алтайский Госуниверситет, Барнаул, 1998.

\item А.А. Карацуба. Основы аналитической теории чисел -- М.: Наука, 1983.

\item E. Landau. Euler and Functionalgleichung der Riemannschen Zeta - Funktion. Biblio-
theca mathematica, (3), Bd.7., H.1,1906, Leipzig, рр.69 -79.

\item A.P. Laurinchikas. Distribution of values of generating Dirichlet series of multiplica-
tive functions. Lit. math.coll., 1982, v.22, 1, pp.101 - 111.

[\item K. Matsumoto. An introduction to the value - distribution theory of zeta - functions. Shiauliai Mathematical Seminar 1(9), 2006, 61 - 83.

\item A. Laurinchikas. Prehistory of the Voronin Universality Theorem. Shiauliai Mathema-
tical Seminar 1(9), 2006, 41-53.

\item Г. Могтномери. Мультипликативная теория чисел.М.Мир,1974.

\item Е. П. Ожигова. Развитие теории чисел в России. Ленинград, Наука, 1972, 360 стр.

\item K. Прахар. Распределение простых чисел. M.: Мир, 1987.

\item B. Riemann. Ueber die Anzahl der Primzahlen unter einer gegebenen Gsцsse.  Monat. der Kцnigl. Preuss. Akad. der Wissen. Zu Berlin. aus der Jahre. 1859(1860), 671-680.

\item У. Рудин. Основы математического анализа. М.: Наука, 1976.

\item Е.К. Титчмарш. Теория дзета-функции Римана. M. : ИЛ, 1953.

\item Е.К. Титчмарш. Теория функций. М.: ГИТТЛ, 1951

\item С.М Воронин, А.А. Карацуба. Дзета-функция Римана. M: физ.-мат.лит, 1994, 376 стр.

\item С.М Воронин. О распределении ненулевых значений дзета- функции Римана // Тр. МИАН -- 1972 -- т. 128, с. 153-175.

\item С.М Воронин. О дифференциальной независимости О -функций // ДАН СССР -- 1973. т. 209, № 6, с.1264 -- 1266.

\item С.М Воронин. О дифференциальной независимости L - функций Дирихле // Аctа Аrith. -1975 , т. ХХVII --с. 493 -- 509.

\item С.М Воронин. Теорема об "универсальности" дзета - функции Римана // Изв. АН СССР сер. мат. -- 1975 --т. 39, № 3 --с. 475 -- 486.

\item С.М Воронин. О нулях дзета - функции квадратичных форм.// Тр. МИАН -- 1976 -- т. 142 -- с. 135 -- 147.

\item С.М Воронин. Аналитические свойства производящих функции Дирихле арифметических объектов: Диссертация на соискание звания д-ра физ.мат.наук МИАН СССР -- М. , 1977 - 90с.

\item С.М Воронин. О нулях некоторых рядов Дирихле, лежащих на критической прямой // Изв. АН СССР сер.мат -1980-т. 44 № 1 -- с. 63 -91.

\item С.М Воронин. О распределении нулей некоторых рядов Дирихле // ТР. МИАН -- 1984 -- т. 163 -- с. 74 -- 77.

\item A. Speiser. Geometrisches zur Riemannschen Zetafunction, Math.Ann.110 (1934), 514 - 521.

\item H. L. Mongomery. The pair correlation of zeroes of the Riemann zeta - function in Analytic Number Theory. (St. Louis), Proc. Sympos. Pure Math. 24. Amer. Math. Soc. Providence, R.I,1973,181-193.

\item J. B. Conrey, A. Chosh, D. Goldston, S. M. Gonek and D. R. Heath - Brown. On the distribution of gaps between zeroes of the zeta - function. Oxford J. Math., Oxford (2), 36, (1985), 43 - 51.

\item А.А. Карацуба. О нижних оценках максимума модуля дзета-функции. Изв. РАН, сер.мат. 2004, 68, 6, 105 - 118.

\item M. Z. Garaev, C. Y. Yildirim. On small distances between ordinates of zeroes of and . ArXiv:math/0610377v2 [math. NT] 19 Mar 2007.

\item M.A. Королев. О больших расстояниях между последовательными нулями дзета-функции Римана. Изв. РАН, сер.мат. 2008, 72, 2, 91 - 104.

\item E. C. Titchmarsh. The theory of the Riemann Zeta Function, 2-nd ed. Revised by H. D. Heath-Brown, Oxford Univ., Press, Oxford, 1986.

\item H. Edwards Riemann's Zeta Function. Academic Press. New York, 1974.

\item E. Bombieri. ``Problems of the Millenium. The Riemann Hypothesis'' CLAY,(2000).

\item B. Conrey. The Riemann Hypothesis. Noticies of AMS, March, 2003, 341-353.

\item P. Sarnak. Problems in Millenium. The Riemann Hypothesis (2004), Princeton Univer. Courant Inc. of Math. Sci.

\item H. Ivaniec and E. Kowalski. ``Analytical Number Theory'' AMS Colloquium Publica-
tions, \textbf{53}, (2004).

\item Jabbarov I. Sh. The Riemann Hypothesis. ArXiv:1006.0381v3, 2010.

\item Jabbarov I. Sh. On the connection between measure and metric in infinite dimensional space. International conference on Differential Equations and Dynamical Systems. Abstracts. Suzdal (Russia) July 2-7, 2010, Moscow, 2010, p. 213-214.

\item Jabbarov I. Sh. On a new measure in infinite dimensional unite sube. ArXiv:1102.3362
v3, 2012.

\item Уолша Дж. Л. Интерполяция и аппроксимация рациональными функциями. ИИЛ, М. 1961.

\item Привалов И. И. Гранические свойства аналитических функций. ГИТТЛ, М.-Л., 1950,

\end{enumerate}

\end{document}